\documentstyle[twoside]{article}
\title{\sc Generalization of a quadratic transformation 
 due to Exton}
\author{\sc Y. S. Kim$^a$,  A. K. Rathie$^b$ and R. B. Paris$^c$\\
\\
${}^a\!$ Department of Mathematics Education, Wonkwang University, Iksan, Korea\\
E-Mail: yspkim@wonkwang.ac.kr\\
${}^b\!$ Department of Mathematics, Central University of Kerala, Kasaragad 671123,\\
Kerala, India\\
E-Mail: akrathie@cukerala.edu.in\\
${}^c\!$ School of Engineering, Computing and Applied Mathematics,\\
 University of Abertay Dundee, Dundee DD1 1HG, UK\\
E-Mail: r.paris@abertay.ac.uk}%
\begin{document}
\def\f#1#2{\mbox{${\textstyle \frac{#1}{#2}}$}}
\newcommand{\fs}{\f{1}{2}}
\newcommand{\g}{\Gamma}
\date{}
\maketitle

\begin{abstract}
Exton [Ganita 54(2003)13-15] obtained numerous new quadratic transformations involving hypergeometric functions of order two and of higher order by applying various known classical summation theorems to a general transformation formula based on the Bailey transformation. We obtain the generalization of one of the Exton quadratic transformations. The results are derived with the help of a generalization of Dixon's summation theorem for the series ${}_3F_2$ obtained earlier by Lavoie {\em et al.} Several interesting known as well as new special cases and limiting cases are also given.

\vspace{0.4cm}

\vspace{0.3cm}

\end{abstract}
\vspace{0.3cm}






\maketitle
\begin{center}
{\bf 1. \  Introduction}
\end{center}
\setcounter{section}{1}
\setcounter{equation}{0}
\renewcommand{\theequation}{\arabic{section}.\arabic{equation}}
The generalized hypergeometric function with $p$ numeratorial and $q$ denominatorial parameters is defined by (see [6, p. 73])
\[{}_pF_q \left[
\begin {array}{c}
    \alpha_1, \ldots, \alpha_p \\ \beta_1, \ldots, \beta_q
\end{array};~z \right]
={}_pF_q \left[
    \alpha_1, \ldots, \alpha_p ; \beta_1, \ldots, \beta_q~;~z
    \right] \]
\begin {equation} \label{1.1}
=\sum_{n=0}^\infty
    \frac{(\alpha_1)_n \ldots (\alpha_p)_n}
         {(\beta_1)_n \ldots (\beta_q)_n}
    \frac{z^n}{n!},
\end{equation}
where $(\alpha)_n$ denotes the Pochhammer symbol (or the shifted factorial, since $(1)_n=n!$) defined for any complex number $\alpha$ by
\[(\alpha)_n=\frac{\Gamma(\alpha+n)}{\Gamma(\alpha)}=
\left\{\begin{array}{ll}\alpha(\alpha+1)\cdots (\alpha+n-1), & n\in {\bf N}=\{1,2,\ldots\}\\
1, & n=0\end{array}\right.\]
When $q=p$ this series converges for all $|x|<\infty$, but when $q=p-1$ convergence occurs when $|x|<1$ (unless the series terminates). 

It should be remarked here that whenever hypergeometric and generalized hypergeometric functions can be summed in terms of  Gamma functions, the results are very important from the application points of view. It should also be noted that summation formulas for ${}_pF_q$ are known for only very restricted arguments and parameters, for example Gauss' two summation theorems, Kummer's summation theorems for the series ${}_2F_1 $, and Dixon's, Watson's, Whipple's and Saalsch$\ddot{\rm u}$tz's summation theorems for the series ${}_3F_2$, and others, play an important role in the theory of hypergeometric and generalized hypergeometric functions. The function ${}_pF_q(z)$ has been extensively studied by many authors such as Slater [7] and Exton [2].

By applying various known summation theorems to a general formula based upon Bailey's transformation theorem given in Slater [7], Exton [3] obtained as a special case numerous new general transformation formulas involving hypergeometric funtions of order two and of higher order. In fact, in our present investigation, we shall be concerned with the following interesting transformation formula
\[\left(\frac{1}{2} + \frac{1}{2}\sqrt{1-x}\right)^{1-2d}
{}_3F_2 \left[ \begin {array}{c}
        2d-1, b, d-\fs \\ 2d-b, d+\fs \end{array};~-\frac{x}{(1+\sqrt{1-x})^2}
        \right]\hspace{2cm}\]
\begin {equation} \label{1.2}
\hspace{6cm}={}_3F_2 \left[ \begin {array}{c}
        d-\fs, d, d-b+\fs \\ 2d-b,d+\fs \end{array};~x\right],
\end{equation}
which is valid for $|x|<1$ and $|x/(1+\sqrt{1-x})^2)|<1$.
Moreover, Exton [3] deduced (1.4) from the following more general transformation formula
\[\left(\frac{1}{2} + \frac{1}{2}\sqrt{1-x}\right)^{1-2d}
{}_{A+1}F_{H+1} \left[ \begin {array}{c}
        (a), d-\fs \\ (h),d+\fs \end{array};~-\frac{xy}{(1+\sqrt{1-x})^2}
        \right]\]
\begin {equation} \label{1.3}
=\sum_{m=0}^\infty\frac{(d-\fs)_m(d)_m}{(2d)_m m!}\,x^m{}_{A+1}F_{H+1} \left[ \begin {array}{c}
        (a), -m \\ (h),2d+m \end{array};~y\right],
\end{equation}
which is valid for $|xy|<1$ and $|xy/(1+\sqrt{1-x})^2|<1$.
Here, the symbol $(h)$ is a convenient contraction for the sequence of parameters $h_1$, $h_2,\ldots, h_H$ and the Pochhammer symbol $(h)_{n}$ is defined above.

 The aim of this paper is to obtain the generalization of (\ref{1.2}) in the form
\begin {equation} \label{1.4}
\left(\frac{1}{2} + \frac{1}{2}\sqrt{1-x}\right)^{1-2d}
{}_3F_2 \left[ \begin {array}{c}
        b, d-\fs, 2d-1-i \\ d+\fs, 2d-b+j \end{array};~-\frac{x}{(1+\sqrt{1-x})^2}\right]
\end{equation}
for $-3\leq i\leq 3$ and $j= 0, {1}, {2}, {3}$.
For this, we shall require the following generalization of Dixon's theorem for the sum of a ${}_3F_2$ of unit argument obtained earlier by Lavoie et al. [4].
\[{}_{3}F_{2}\left[\begin{array}{c} a,\qquad\qquad
b,\qquad\qquad c \\
       1+a-b+i, 1+a-c+i+j \end{array}; 1\right]\hspace{5cm} \]     
\[={2^{-2c+i+j}\Gamma(1+a-b+i)\Gamma(1+a-c+i+j)\Gamma(b-{1\over2}i-{1\over2}|i|)\Gamma(c-{1\over2}(i+j+|i+j|))
               \over\Gamma(a-2c+i+j+1)\Gamma(a-b-c+i+j+1)\Gamma(b)\Gamma(c)}\]
\[\times \Biggl\{A_{i,j}{\Gamma(\frac12 a-c+{1\over2}+[\frac{i+j+1}{2}])\Gamma(\frac12 a-b-c+1+i+[\frac{j+1}{2}])
               \over\Gamma(\frac12 a+\frac12)\Gamma(\frac12 a-b+1+[\frac12 i])}\]
\begin{equation}\label{1.5}
        +B_{i,j}{\Gamma(\frac12 a-c+1+[\frac{i+j}{2}])\Gamma(\frac12 a-b-c+\frac32+i+[\frac{j}{2}])
               \over\Gamma(\frac12 a)\Gamma(\frac12 a-b+\frac12+[\frac{i+1}{2}])}\Biggr\}
\end{equation}
provided $\Re (a-2b-2c)>-2-2i-j$ with $-3\leq i\leq 3$  and $j= 0,{1}, {2}, {3}$.
Here and in what follows, $[x]$ is the greatest integer less than or equal to $x$ and $|x|$ denotes the usual absolute value of $x$.
The coefficients $A_{i,j}$ and $B_{i,j}$ are given in Tables 1 and 2.
\begin{table}[th]
\scriptsize
\caption{Values of the coefficients $A_{i,j}$}
\begin{center}
\begin{tabular}{|c|c|c|c|c|} 
$i\backslash j$ & 0 & 1 & 2 & 3
\\  \hline
3 & $\begin{array}{c}
               5a-b^2+(a+1)^2 \\ -(2a-b+1)
               (b+c) \end{array}$ & {\rm ---} & {\rm ---} & {\rm ---} \\ \hline
2 & $\begin{array}{c}
              \frac12\, (a-1)(a-4) \\ - (b^2-5a+1)\\ - (a-b+1)
               (b+c) \end{array}$  &
           $\begin{array}{c}
               (b-1)(b-2) \\- (a-b+1) \\ \times (a-b
               -c+3) \end{array}$ &
$\begin{array}{c}
 \fs\,(a-c+2)
 (a-2b-c+5) \\
 \times \{ (a-c+2)
 (a-2b+2)\\
 -a(c-3)\} \\
 -(b-1)(b-2)
 (c-2)(c-3) \end{array}$ & {\rm ---}
 \\ \hline
1 & 1  & $c-a-1$ & $\begin{array}{c}
                              a(a-1)\\+ (b+c-3) (c-2a-1) \end{array}$ & {\rm ---} \\ \hline
 0 & 1  & $-1$ & $
 \begin{array}{c}
                 \fs\,\{(a-b-c+1)^2 \\ +(c-1)(c-3)-b^2+a\} \end{array}$
& $\begin{array}{c}
                         c(a-b-c+4) \\- (a+1)(a+2)\\-(a-1)(b-1)+3ab\end{array}$ \\ \hline
$-1$ &1 &1 &$b+c-1$
    &$\begin {array}{c} (c-1)(c-2)\\-b(a-c+1)
     \end {array}$\\ \hline
$-2$ & $\begin {array}{c}
    \fs (a-1) (a-2b-2) \\ -c(a-b-1)
    \end {array}$
    &$a-b-1$
    &$\begin {array}{c} \fs (a-1)(a-2b-2c)
    \\ +b(b+c) \end {array}$
    &$\begin {array}{c} (a-b-1)(c-1) \\
     -b(b+1) \end {array}$ \\ \hline
$-3$ &$\begin {array}{c}(a-1) \\\times(a-2b-2c-4)\\
     +bc \end {array}$
    &$\begin {array}{c} (a-b-2)\\\times (a-c-1) \\
    -ac \end {array}$
    &$\begin {array}{c} (a-b-1) \\
    \times (a-b-2c-2) \\
    -bc \end {array}$
    &$\begin {array}{c} b(b+1) \\+(a-1)(a-b) \\
    -c(2a-b-2) \end {array}$ \\ \hline
\end{tabular}
\end{center}
\end{table}
\begin{table}[h]
\scriptsize
\caption{Values of the coefficients $B_{i,j}$}
\begin{center}
\begin{tabular}{|c|c|c|c|c|} 
        $i\backslash j$ & 0 & 1 & 2 & 3 \\ \hline
 3 &
 $\begin{array}{c}
 -a+3b^2
 -(a+3)^2 \\ +(2a-3b+5)(b+c) \end{array}$ & {\rm ---} & {\rm ---} & {\rm ---} \\ \hline
2 & -2 & $\begin{array}{c}
                    (a-b-2c+5) \\ \times (a-b-c+3) \\ -(b-1)(b-2) \end{array}$
& $\begin{array}{c}
-2(a-c+2)\\
\times (a-2b-c+5)
\end{array}$& {\rm ---} \\ \hline
1 & -1 & $a-2b-c+3$ & $\begin{array}{c}
           (b-1) (b-c+1) \\-(a-b-c+2)\\
\times (a-b-c+3) \end{array}$ & {\rm ---} \\ \hline
       0 & 0  & 1 & -2
& $\begin{array}{c}
 (a+2)(a+4)-b(2a+5) \\- 3c(a-b-c+4) +3 \end{array}$ \\ \hline
$-1$& 1 &1 &$-(b-c+1)$
  &$\begin {array}{c} (c-1(c-2) \\
  +b (a-2b-c+1) \end {array}$ \\ \hline
$-2$ &2 &$a-b-2c-1$ & 2
    &$\begin {array}{c} b(a-2c+2) \\
    -(b-c+1) \\ \times (a-b-2c+1) \end {array}$ \\
    \hline
$-3$ &$\begin {array}{c} (a-2)\\ \times(a-2b -2c-3) \\
    +3bc \end {array}$
    &$\begin{array}{c} (a-b-2) \\
    \times (a-2b \\-2c-3)+bc  \\ \end{array}$
    &$\begin {array}{c} (a-b-2) \\ \times (a-b-2c-1)\\+bc \\
    \end {array}$
    &$\begin {array}{c} (a-1)(a-2) \\ -3b(a-b-2) \\
    -c(2a-3b-4) \end {array}$ \\ \hline
\end{tabular}
\end{center}
\end{table}
\noindent Also, if $f_{i,j}$ denotes the ${}_3F_2(1)$ series on the left-hand side of (\ref{1.5}), the natural symmetry
\[f_{i,j}(a, b, c) = f_{i+j,-j}(a, c, b)\]
makes it possible to extend the result to $j= {-1}, {-2}, {-3}$.

Several interesting cases, including Exton's result, are then deduced as special cases of our main findings. In addition to this, certain known results obtained recently by Pog\'any and Rathie [5] are also obtained as a limiting case of our main findings. The results derived in this paper are simple, interesting, easily established and may be useful.
\vspace{0.6cm}

\begin{center}
{\bf 2. \  Extension of Exton's quadratic transformation}
\end{center}
\setcounter{section}{2}
\setcounter{equation}{0}
\renewcommand{\theequation}{\arabic{section}.\arabic{equation}}
Here we establish a natural extension of the Exton transformation (\ref{1.2}) given by the following theorem.
\newtheorem{theorem}{Theorem}
\begin {theorem}
The following identities hold true in the domain ${\cal D}$ defined by the connected subset ${\cal D} =\left\{x\in {\bf C}\,\,\big|\,\,|x|<1,  |x/(1+\sqrt{1-x})^2|<1\right\}$:
\[\left(\frac{1}{2} + \frac{1}{2}\sqrt{1-x}\right)^{1-2d}
{}_3F_2 \left[ \begin {array}{c}
        b, d-\fs, 2d-1-i \\ d+\fs, 2d-b+j \end{array};~-\frac{x}{(1+\sqrt{1-x})^2}
        \right]\]
 \[={2^{i}(-1)^{{1\over2}(i+|i|)}\Gamma(d)\Gamma(d+\frac12)\Gamma(b-{1\over2}(i+j+|i+j|))
               \over\Gamma(b)\Gamma(d-b+{\frac12}j)\Gamma(d-b+\fs j+\fs)}\\
   \times \sum_{n=0}^{\infty} {n!\over (n+{1\over2}i+{1\over2}|i|)!}{(d)_n(d-\frac 12)_n \over(2d-b+j)_n}{x^n \over n!}\]
 \[  \times \Biggl\{A_{i,j}{\Gamma(d-b-\frac{i}{2}+[\frac{i+j+1}{2}])\Gamma(d-b+\frac{i}{2} +\frac12+[\frac{j+1}{2}])
               \over\Gamma(d-\frac{i}{2})\Gamma(d+\frac 12-\frac{i}{2}+[\frac{i}{2}])}{(d-b+\frac{i}{2}+\frac12+[\frac{j+1}{2}])_n \over (d+\frac 12-\frac{i}{2}+[\frac{i}{2}])_n}\]
 \begin {equation} \label{2.1}
        +B_{i,j}{\Gamma(d-b+\fs-\frac{i}{2}+[\frac{i+j}{2}])\Gamma(d-b+\frac{i}{2} +1+[\frac{j}{2}])
               \over\Gamma(d-\frac{i}{2}-\frac12)\Gamma(d-\frac{i}{2}+[\frac{i+1}{2}])}{(d-b+\frac{i}{2}+1+[\frac{j}{2}])_n \over (d-\frac{i}{2}+[\frac{i+1}{2}])_n}\Biggr\}
\end{equation}
for $-3\leq i\leq 3$ and $j= 0,{1}, {2}, {3}$.
As usual  $[x]$ denotes the greatest integer less than or equal to $x$ and its modulus is denoted by $|x|$.
The coefficients $A_{i,j}$ and $B_{i,j}$ can be obtained from the values of $A_{i,j}$ and $B_{i,j}$ in Tables 1 and 2 by changing $a$ to $2d-1-i$, $b$ to $-n$ and $c$ to $b$, respectively.
\end{theorem}
\bigskip

\noindent{\em Proof}\ \ \ 
 We first derive Exton's result (\ref{1.3}) in an alternative way. Let $\mathcal{S}$ denote the left-hand side of (\ref{1.3}) and express ${}_{A+1}F_{H+1}$ as a series so that
\[S = \sum_{n=0}^{\infty}{(-1)^n((a))_n (d-\fs)_n x^{n}y^{n}\over ((h))_n (d+\fs)_n 2^{2n}n!}{\left(\frac{1}{2} + \frac{1}{2}\sqrt{1-x}\right)^{1-2(d+n)}}.
\]
Use of the well-known result [8, p.~34]
\[\left(\frac{1}{2} + \frac{1}{2}\sqrt{1-x}\right)^{1-2a}
={}_2F_1 \left[ \begin {array}{c}
        a-\fs,\,\, a \\ 2a \end{array};~x\right],\]
then enables $S$ to be rewritten in the form
\[S = \sum_{n=0}^{\infty}{(-1)^n((a))_n (d-\fs)_n x^{n}y^{n}\over ((h))_n (d+\fs)_n 2^{2n}n!}\,{}_2F_1 \left[ \begin {array}{c}
       d+ n-\fs,\,\, d+n \\ 2d+2n \end{array};~x\right].\]

Expressing ${}_2F_1$ as a series, we obtain
\[
S = \sum_{n=0}^{\infty}\sum_{m=0}^{\infty}{(-1)^n((a))_n (d-\fs)_n (d+n-\fs)_m (d+n)_m x^{n+m}y^{n}\over ((h))_n (d+\fs)_n (2d+2n)_m 2^{2n}n!m!}.\]
 Changing $m$ to $m-n$ and using the following identities
[6, p.~57, Eq.~(8)]
\[\sum_{n=0}^{\infty}\sum_{k=0}^{\infty}A(k, n) =
\sum_{n=0}^{\infty}\sum_{k=0}^{n}A(k, n-k)\]
and
\[(\alpha+k)_{n-k} = \frac{(\alpha)_n}{(\alpha)_k},\qquad
(n-k)!=\frac{(-1)^k~ n!}{(-n)_k},\]
we find after some simplification that
\[
S = \sum_{m=0}^{\infty}{(d-\fs)_m (d)_m x^{m}\over  (2d)_m m!}\sum_{n=0}^{m}{((a))_n (-m)_n y^n \over ((h))_n (2d+m)_n n!}.\]
Finally, summing the inner series as a hypergeometric series, we easily arrive at the right-hand side of (\ref{1.3}). This completes our proof of (\ref{1.3}).

Now we are ready to derive our main result (2.1). For this, if we put $A = 2$, $H = 1$, $a_{1} = 2d-1-i$, $a_{2} = b$, $h_{1} = 2d-b-j$ and $y = 1$ in (\ref{1.3}), we obtain
\[
\left(\frac{1}{2} + \frac{1}{2}\sqrt{1-x}\right)^{1-2d}
{}_3F_2 \left[ \begin {array}{c}
        2d-1-i, b, d-\fs \\ 2d-b+j,d+\fs \end{array};~-\frac{x}{(1+\sqrt{1-x})^2}
        \right]\]
\begin{equation}\label{2.2}
=\sum_{n=0}^{\infty}{(d-\fs)_n (d)_n x^{n}\over  (2d)_n n!}\,{}_3F_2 \left[ \begin {array}{c}
        2d-1-i,\,\, b,\,\, -n \\ 2d-b+j,\,\,2d+n \end{array};~1
         \right].
\end{equation}
It is easy now to see that the ${}_3F_2$ on the right-hand side of (\ref{2.2}) can be evaluated with the help of the generalized Dixon summation theorem (\ref{1.5}) by taking $a$ by $2d-1-i$, $b$ by $-n$ and $c$ by $b$, and  after a little simplification, we easily arrive at the right-hand side of (\ref{2.1}). This completes the proof of (\ref{2.1}).
\hfill $\Box$
\vspace{0.6cm}

\begin{center}
{\bf 3. \  Special cases}
\end{center}
\setcounter{section}{3}
\setcounter{equation}{0}
\renewcommand{\theequation}{\arabic{section}.\arabic{equation}}
By assigning values to $i$ and $j$ in our main result (\ref{2.1}), we can obtain a large number of interesting and useful results. However, we shall mention here only a few of them. All these transformations hold in a domain ${\cal D}$ defined by the connected subset ${\cal D} =\{x\in {\bf C}\,\,|\,\,|x|<1,  |x/(1+\sqrt{1-x})^2|<1\}$.

For $i=0$ and $j=0$ in (\ref{2.1}), we obtain 
\[\left(\frac{1}{2} + \frac{1}{2}\sqrt{1-x}\right)^{1-2d}
{}_3F_2 \left[ \begin {array}{c}
        2d-1,\, b,\, d-\fs \\ 2d-b,\,d+\fs \end{array};~-\frac{x}{(1+\sqrt{1-x})^2}\right]\hspace{3cm}\]
\begin{equation} \label{3.1}
\hspace{5cm}={}_3F_2 \left[ \begin {array}{c}
        d-\fs,\, d,\, d-b+\fs \\ 2d-b,\,d+\fs \end{array};~x
         \right],
\end{equation}
which is the result stated in (\ref{1.2}).

 For $i=0$ and $j=1$ in (\ref{2.1}), we obtain
\[\left(\frac{1}{2} + \frac{1}{2}\sqrt{1-x}\right)^{1-2d}
{}_3F_2 \left[ \begin {array}{c}
        2d-1,\, b,\, d-\fs \\ 2d-b+1,\,d+\fs \end{array};~-\frac{x}{(1+\sqrt{1-x})^2}
        \right]\]
\[={(2d-2b+1)\over 2(1-b)}\,{}_3F_2 \left[ \begin {array}{c}
        d-\fs,\, d,\, d-b+\f{3}{2} \\ 2d-b+1,\,d+\fs \end{array};~x
         \right] \]
\begin{equation} \label{3.2}
 -{(2d-1)\over 2(1-b)}\,{}_2F_1 \left[ \begin {array}{c}
        d-\fs,\,\, d-b+1 \\ 2d-b+1 \end{array};~x
         \right].\hspace{1cm}
\end{equation}

For $i=1$ and $j=0$ in (2.1), we obtain
\[
\left(\frac{1}{2} + \frac{1}{2}\sqrt{1-x}\right)^{1-2d}
{}_3F_2 \left[ \begin {array}{c}
        2d-2,\, b,\, d-\fs \\ 2d-b,\,d+\fs \end{array};~-\frac{x}{(1+\sqrt{1-x})^2}
        \right]\]
\[={(2d-1)(d-b)\over (1-b)}\,{}_3F_2 \left[ \begin {array}{c}
        d-\fs,\, d-b+1,\, 1 \\ 2d-b,\, 2 \end{array};~x
         \right] \]
\begin{equation} \label{3.3}
-{(d-1)(2d-2b+1)\over (1-b)}\,{}_4F_3 \left[ \begin {array}{c}
        d,\,d-\fs,\, d-b+\f{3}{2},\,1 \\ 2d-b,\, d+\fs,\,2 \end{array};~x
         \right].
\end{equation}

 For $i=1$ and $j=1$ in (2.1), we obtain
\[\left(\frac{1}{2} + \frac{1}{2}\sqrt{1-x}\right)^{1-2d}
{}_3F_2 \left[ \begin {array}{c}
        2d-2,\, b,\, d-\fs \\ 2d-b+1,\,d+\fs \end{array};~-\frac{x}{(1+\sqrt{1-x})^2}
        \right]\]
\[={(2d-1)(d-b+1)(2d-b-1)\over (b-1)(b-2)}\,{}_3F_2 \left[ \begin {array}{c}
        d-\fs,\, d-b+2,\, 1 \\ 2d-b+1,\,2 \end{array};~x
         \right] \]
\begin{equation} \label{3.4}
-{(d-1)(2d-b+1)(2d-2b+1)\over (b-1)(b-2)}\,{}_5F_4 \left[ \begin {array}{c}
        d-\fs,\,d,\, d-b+\f{3}{2},\,d-\fs b+\f{3}{2},\,1 \\ 2d-b+1,\,d+\fs,\,d-\fs b+\fs,\,2 \end{array};~x
         \right].
\end{equation}

For $i=-1$ and $j=0$ in (2.1), we obtain
\[\left(\frac{1}{2} + \frac{1}{2}\sqrt{1-x}\right)^{1-2d}
{}_3F_2 \left[ \begin {array}{c}
        2d,\, b,\, d-\fs \\ 2d-b,\,d+\fs \end{array};~-\frac{x}{(1+\sqrt{1-x})^2}
        \right]\]
\begin{equation} \label{3.5}
=\frac{1}{2}\,{}_2F_1 \left[ \begin {array}{c}
        d-\fs,\, d-b \\ 2d-b \end{array};~x
         \right]
 +\frac{1}{2}\,{}_3F_2 \left[ \begin {array}{c}
       \,d,\, d-\fs,\, d-b+\fs \\ 2d-b,\,d+\fs \end{array};~x
         \right].
\end{equation}

For $i=-1$ and $j=1$ in (2.1), we obtain
\[\left(\frac{1}{2} + \frac{1}{2}\sqrt{1-x}\right)^{1-2d}
{}_3F_2 \left[ \begin {array}{c}
        2d, b, d-\fs \\ 2d-b+1,\,\,d+\fs \end{array};~-\frac{x}{(1+\sqrt{1-x})^2}
        \right]\]
\begin{equation} \label{3.6}
=\frac{1}{2}\,{}_2F_1 \left[ \begin {array}{c}
        d-\fs,\, d-b+1 \\ 2d-b+1 \end{array};~x
         \right]
 +\frac{1}{2}\,{}_3F_2 \left[ \begin {array}{c}
       \,d,\, d-\fs,\, d-b+\fs \\ 2d-b+1,\,d+\fs \end{array};~x
         \right].
\end{equation}

For $i=-2$ and $j=1$ in (2.1), we obtain
\[\left(\frac{1}{2} + \frac{1}{2}\sqrt{1-x}\right)^{1-2d}
{}_3F_2 \left[ \begin {array}{c}
        2d+1, b, d-\fs \\ 2d-b+1,\,\,d+\fs \end{array};~-\frac{x}{(1+\sqrt{1-x})^2}
        \right]\]
\begin{equation} \label{3.7}
=\frac{1}{2}\,{}_4F_3 \left[ \begin {array}{c}
        d,\,d-\fs,\,2d+1,\, d-b+\fs \\ 2d,\, d+\fs,\,2d-b+1 \end{array};~x
         \right]
+\frac{1}{2}\,{}_3F_2 \left[\begin{array}{c}
       \,d-b,\, d-\fs,\, 2d-2b+1 \\ 2d-b+1,\,2d-2b \end{array};~x
         \right].
\end{equation}
Similarly other results can also be obtained.
\vspace{0.6cm}

\begin{center}
{\bf 4. \  Limiting cases}
\end{center}
\setcounter{section}{4}
\setcounter{equation}{0}
\renewcommand{\theequation}{\arabic{section}.\arabic{equation}}
Here we mention some of the interesting limiting cases of our results. All these transformations hold in the domain $\mathcal{D}$ defined by the connected subset ${\cal D} =\{x\in {\bf C}\,\,\big|\,\,|x|<1,  |x/(1+\sqrt{1-x})^2|<1\}$.

If we let $b \rightarrow \infty$ in (\ref{3.1}) or (\ref{3.2}), we obtain the following result:
\[\left(\frac{1}{2} + \frac{1}{2}\sqrt{1-x}\right)^{1-2d}
{}_2F_1 \left[ \begin {array}{c}
        2d-1,\,\, d-\fs \\\,d+\fs \end{array}~;~\frac{x}{(1+\sqrt{1-x})^2}
        \right]\]
\begin {equation} \label{4.1}
={}_2F_1 \left[ \begin {array}{c}
        d-\fs,\,\, d \\ \,d+\fs \end{array};~x
         \right].
\end{equation}

If we let $b \rightarrow \infty$ in (\ref{3.3}) or (\ref{3.4}), we obtain the following result:
\[\left(\frac{1}{2} + \frac{1}{2}\sqrt{1-x}\right)^{1-2d}
{}_2F_1 \left[ \begin {array}{c}
        2d-2, \,\, d-\fs \\ \,d+\fs \end{array};~\frac{x}{(1+\sqrt{1-x})^2}
        \right]\]
\begin {equation} \label{4.2}
=(2d-1)\,{}_2F_1 \left[ \begin {array}{c}
        d-\fs, 1 \\ 2 \end{array};~x
         \right]
 -2(d-1)\,{}_3F_2 \left[ \begin {array}{c}
        d-\fs,\,d,\, 1 \\ d+\fs,\,2 \end{array};~x
         \right].
\end{equation}

If we let $b \rightarrow \infty$ in (\ref{3.5}) or (\ref{3.6}), we obtain the following result:
\[\left(\frac{1}{2} + \frac{1}{2}\sqrt{1-x}\right)^{1-2d}
{}_2F_1 \left[ \begin {array}{c}
        2d, \,\, d-\fs \\ \,d+\fs \end{array};~\frac{x}{(1+\sqrt{1-x})^2}
        \right]\]
\begin {equation} \label{4.3}
=\frac{1}{2}\,{}_1F_0 \left[ \begin {array}{c}
        d-\fs \\ - \end{array};~x
         \right]
 +\frac{1}{2}\,{}_2F_1 \left[ \begin {array}{c}
        d-\fs,\,d\\ d+\fs \end{array};~x
         \right].
\end{equation}

If we let $b \rightarrow \infty$ in (\ref{3.7}), we obtain the following result:
\[\left(\frac{1}{2} + \frac{1}{2}\sqrt{1-x}\right)^{1-2d}
{}_2F_1 \left[ \begin {array}{c}
        2d+1, \,\, d-\fs \\ \,d+\fs \end{array};~\frac{x}{(1+\sqrt{1-x})^2}
        \right]\]
\begin {equation} \label{4.4}
=\frac{1}{2}\,{}_1F_0 \left[ \begin {array}{c}
        d-\fs \\ - \end{array};~x
         \right]
 +\frac{1}{2}\,{}_3F_2 \left[ \begin {array}{c}
        d-\fs,\,d,\,2d+1\\ d+\fs,\,2d \end{array};~x
         \right].
\end{equation}

We conclude this section with a remark that the result (\ref{4.1}) was obtained by Choi and Rathie [1] whereas the results (\ref{4.2})--(\ref{4.4}) were obtained by Pog\'any and Rathie [5] by using a generalization of Kummer's summation theorem. For a remark on the Exton result [3], see the paper by Choi and Rathie [1].
\bigskip

\centerline{\textbf{Acknowledgement}}\medskip
The first author acknowledges the financial support of  Wonkwang University in 2014.

\end{document}